\documentclass[a4paper]{article}

\usepackage{amsmath}
\usepackage{amsfonts}
\usepackage{graphicx}
\usepackage[margin=3.5cm]{geometry}

\newtheorem{thm}{Theorem}

\begin{document}

\title{A Concise, Elementary Proof of \\Arzel\`a's Bounded Convergence Theorem}
\author{Nadish de Silva}
\date{}
\maketitle

\begin{abstract}
Arzel\`a's bounded convergence theorem (1885) states that if a sequence of Riemann integrable functions on a closed interval is uniformly bounded and has an integrable pointwise limit, then the sequence of their integrals tends to the integral of the limit.  It is a trivial consequence of measure theory.  However, denying oneself this machinery transforms this intuitive result into a surprisingly difficult problem; indeed, the proofs first offered by Arzel\`a and Hausdorff were long, difficult, and contained gaps.  In addition, the proof is omitted from most introductory analysis texts despite the result's naturality and applicability.  Here, we present a novel argument suitable for consumption by freshmen.
\end{abstract}

\section{Introduction}

The proofs of Arzel\`a's bounded convergence theorem first offered by Arzel\`a and Hausdorff were long, difficult, and contained gaps.  Indeed, most introductory analysis texts record the result but do not supply a proof.  For an excellent history of the problem and a cogent explanation of the hypotheses' necessity see \cite{gordon}.  Many other proofs that do not make use of measure theory have been furnished by mathematicians such as Banach and Reisz, using techniques ranging from set theory to functional analysis---see \cite{lux} for more about these proofs.  Our goal is to establish the following theorem using only Riemann's definition of integration and basic principles of elementary analysis:

\begin{thm}\label{premain}
\emph{(Arzel\`a, 1885)} Let $f_n:[a,b]\rightarrow\mathbb{R}$ be a sequence of Riemann integrable functions with the pointwise limit $f:[a,b]\rightarrow\ \mathbb{R}$ which is also integrable.  Suppose further that there is a constant $C$ bounding the $f_n$ uniformly: $|f_n(x)| < C$ for every $x \in [a,b]$ and $n \in \mathbb{N}$.  Then $\int_a^b f = \lim_{n\to\infty} \int_a^b f_n$.

\end{thm}

\noindent By considering the related sequence of functions $\{|f_n - f|\}$ and rescaling, one sees that we may assume that both the domain and codomain of our functions are the unit interval and that the limit function is zero:

\begin{thm}\label{main}
Let $f_n:[0,1]\rightarrow[0,1]$ be a sequence of integrable functions that tend pointwise to zero.  Then $\lim_{n\to\infty} \int_0^1 f_n=0$.

\end{thm}

\section{Proof}

\noindent We will prove the contrapositive.  By passing to a subsequence if necessary, we may assume without loss of generality that the integrals $\int_0^1 f_n$ are all above a fixed bound.  We will find a point $x \in [0,1]$ such that infinitely many $f_n(x)$ are also above a fixed bound.  We conclude that neither the subsequence nor the original sequence converges pointwise to zero.

So, suppose there is an $\epsilon$ such that $\int f_n > 2\epsilon$ for every $n \in \mathbb{N}$. We can associate to each $f_n$ a finite union $U_n$ of open intervals whose total length is at least $\epsilon$ such that $f_n$ only takes on values above $\epsilon$ on $U_n$.  
To see this, note that the upper/lower sum definition of the Riemann integral tells us that we may inscribe a finite number of rectangles under the graph of $f_n$ whose combined area is at least  $2\epsilon$.  
Since the short rectangles (those whose height is less than $\epsilon$)  have a total area of at most $\epsilon$, the total area of the tall rectangles is at least $\epsilon$.
Taking $U_n$ to be the open intervals that constitute the bases of these tall rectangles, we conclude that the length of $U_n$ is at least $\epsilon$ since the height of any rectangle is at most one.  It is clear that $f_n$ is above $\epsilon$ on $U_n$.  

\begin{figure}[h!]
  \centering
\includegraphics[width=2in]{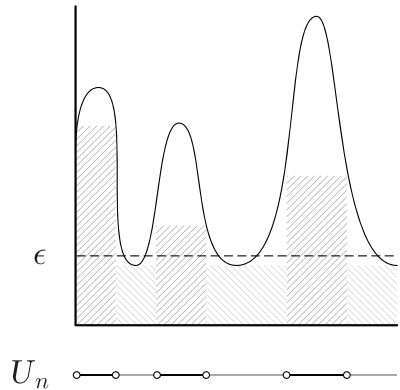} 
  \caption{To each $f_n$ we associate a set $U_n$ which is the union of finitely many open intervals.  The length of $U_n$ is greater than $\epsilon$, as are the values taken by $f_n$ on $U_n$.}
\end{figure}

Finding a point common to infinitely many $U_n$ gives us the desired point $x$ of nonconvergence to zero.  Equivalently, we will show that if $V_n = \bigcup_{k=n}^\infty U_k$ , then the intersection of all $V_n$ contains a point.

\begin{thm}\label{secondary}
Let $\{V_n\}$ be a nested, decreasing sequence of open sets in $[0,1]$, each of which contains a finite union of open intervals whose total length is above some fixed bound $\epsilon$.  Then $\bigcap_{n=1}^\infty V_n$ is nonempty.
\end{thm}

Measure theory immediately tells us that the intersection is not merely nonempty but, indeed, of positive measure.  That we must settle for such a weak conclusion is a testament to the theory's power.  To establish this theorem by elementary means, we introduce a natural tree structure on the sequence $\{V_n\}$.  As $V_n$ is open, it decomposes uniquely as a countable disjoint union of open intervals.  Each such interval is a node at depth $n$; we speak of the properties of each interchangeably.  The children of a node $N$, an interval of $V_n$, are the intervals of $V_{n+1}$ that are subintervals of $N$.  We will prune this tree until the sequence of open sets which remains can be easily seen to have a nonempty intersection.

\begin{figure}[h!]
  \centering
\includegraphics[width=2in]{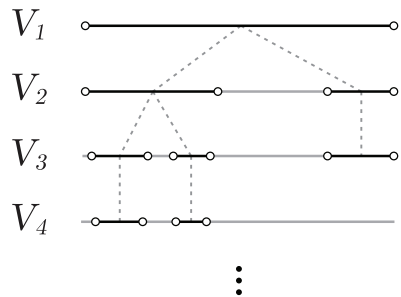} 
  \caption{A nested sequence of nonempty open sets can be thought of as a tree with infinite depth.}
\end{figure}

We may assume without loss of generality that each $V_n$ is a finite union of intervals.  To see this, enumerate the intervals that constitute $V_n$ and notice that the total length of the first $k$ intervals is a monotone sequence bounded above by one.  We may thus truncate the sequence of intervals such that the length of the resulting finite union of intervals is arbitrarily close to the limit length.  Truncate each $V_n$ to within ${\epsilon}/{2^{n+1}}$ of this limit.  By discarding the descendants of the truncated tails we ensure that the resulting sequence of open sets is still nested.  The lengths of the finite unions of intervals that remain at each depth are all at least ${\epsilon}/{2}$.

We may further assume that there are no \emph{terminating nodes} (nodes whose descendants end at a finite depth).  They may be safely discarded as the nonterminating nodes of depth $n$ wholly contain $V_{n'}$ for sufficiently large $n'$.  So, all intervals can be assumed to have descendants at every deeper level.

The case where there is an infinite path down the tree such that the lengths of the intervals traversed by the path stays above a fixed number is easy.  The midpoints of the intervals traversed have a cluster point; it is in the intersection of all the $V_n$s.  We may thus assume that as we travel down any path, the lengths of the intervals traversed vanish.

One last reduction: all the nodes \emph{split} (have at least two descendants at sufficiently large depths).  Each $V_n$ is the disjoint union of $S_n$, the splitting nodes, and $O_n$, the nodes with only one descendant at each further depth.  The descendants of $O_n$ at depth $n'$ vanish in length as $n'$ increases.  This means the descendants of $S_n$, and thus $S_n$ itself, are longer than some fixed uniform bound.  We can discard $O_n$ at every level.

Now, every node is nonterminating and splits.  Note that any node $N$ splits into two descendants and that one of these subintervals splits, implying that $N$ has three descendants at some depth.  The midpoints of the outermost of these nodes define a closed subinterval of $N$ that wholly contains the middle descendant.  Starting with any initial node and inductively applying this construction to each such yielded middle node gives a nested sequence of closed sets whose nonempty intersection is contained in $\bigcap_{n=1}^\infty V_n$.

\begin{figure}[h!]
  \centering
\includegraphics[width=4in]{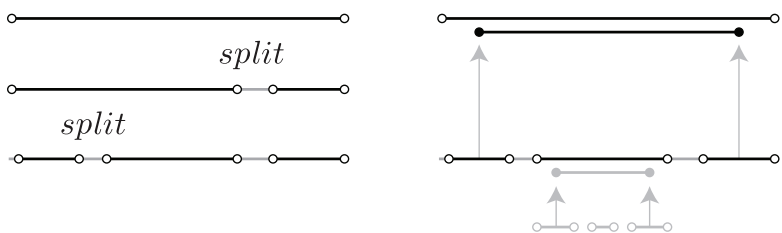} 
  \caption{In a tree where every node is nonterminating and splits, every interval contains a closed set which in turn contains a descendant node.}
\end{figure}

\noindent {\bf Acknowledgements.} I would like to thank Andres del Junco for interesting me in this problem during a real analysis class.  In addition, I am grateful to George Elliott and Dror Bar-Natan for their guidance and encouragement and to Alex Bloemendal, Larry Guth, and the referees for their helpful comments.  I would also like to acknowledge the support of the the Department of Mathematics at the University of Toronto.

\bigskip

\noindent\textit{Department of Mathematics, 
University of Toronto, Toronto, CA\\
Oxford University Computing Laboratory, University of Oxford, Oxford, UK\\
nadish.desilva@utoronto.ca}

\end{document}